\documentclass[a4paper]{article}
\usepackage[english]{babel}
\usepackage[T1]{fontenc}
\usepackage[utf8]{inputenc}
\usepackage{amsmath,amssymb}
\usepackage{graphicx}
\usepackage[colorinlistoftodos]{todonotes}
\usepackage{hyperref}
\graphicspath{ {./images/} }

\newtheorem{thm}{Theorem}
\newtheorem{prop}{Proposition}

\title{Random Walk in Random Environment: A short introduction}
\author{Mouad El Bouchattaoui \thanks{Currently affiliated with \href{mailto:mouad.ebouchattaoui@centralesupelec.fr}{Paris-Saclay University, CentraleSupélec}.}\\Ecole Centrale Casablanca 
\\\href{mailto:mouad.elbouchattaoui@centrale-casablanca.ma}{mouad.elbouchattaoui@centrale-casablanca.ma}}

\date{05 January 2019}

\begin{document}
\maketitle
\tableofcontents
% \newpage
\begin{abstract}
This is a report of a scientific project carried out at Ecole Centrale Casablanca in 2019. This work is an entry into the world of Random Walk in a Random Environment (RWRE). We will discuss some intuitive concepts such as recurrence, location, and random potential. We will thereafter focus on providing an "electrical" understanding of this type of walk and finally explore an application of RWRE to graphene conduction. 
\end{abstract}

\section{Introduction}
\label{sec:introduction}
Random Walk in a Random Environment (RWRE) is a mathematical subject that is mainly related to probability theory. It is a subject that comes from physics and biology. For example, in physics, it is a powerful tool for modeling the diffusion of particles in a material with irregularities and defects. In biology, it can be used to model the replication of DNA chains.

Intuitively, the concept of RWRE is well understood in random media theory, where the behavior of a particle is highly dependent on the local properties of the diffusion medium. However, understanding these local properties in each element of the material is difficult because of the defects of the material that cannot be located or classified one by one as the "walker" travels through the material. This lack of environmental information is modeled by the idea of the random environment.

Therefore, the random environment is obtained by choosing the local characteristics of the movement randomly and following a probability law. But before the walk begins to evolve over time, the properties of the environment are randomly selected, i.e., a random environment is randomly set. During the evolution, the transition from one site to another is made by a transition probability.

\section{RWRE: A Motivation}
\label{sec:A Motivation}

The model of RWRE was first introduced in 1967 by Chernov \cite{Chernov} to model DNA chain replication. This is a one-dimensional model where the walker can only jump to the nearest nearby sites in each step.

\begin{figure}[!htbp]
    \centering
    \includegraphics{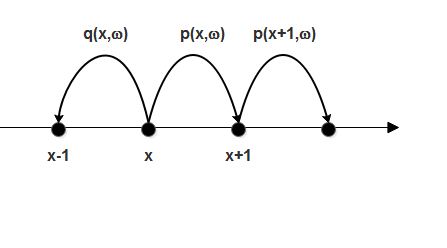}
    \caption{Chernov Model}
    \label{1D_walk}
\end{figure}

In the Chernov model, we choose a family of i.i.d. random variables \begin{math} (p(x,\omega))_{x \in \mathbb{Z}} \end{math} which represents, being in site \( x \), the probability of jumping to the site \( x+1 \) (cf. Fig \ref{1D_walk}). The \begin{math} q(x,\omega) \end{math} is the probability of going to \( x-1 \).
We have thus the relation: \begin{math} q(x,\omega)= 1-p(x,\omega) \end{math}.
By fixing an environment \begin{math} \omega \end{math}, the trajectory of the random walk in each step \( n \), \begin{math} X_{n} \end{math}, defines a Markov chain. In the following section, we will define in detail the RWRE on \(\mathbb{Z}\).

\section{Definition of RWRE on \texorpdfstring{$\mathbb{Z}$}{the set of integers}}

\subsection{Quenched Law and Annealed Law}
We know now that an environment is chosen randomly and remains fixed during the temporal evolution of the walk.
We denote by \(\Omega\) the collection of possible environments that will be equipped with a probability law noted \(\mathbb{P}\). By fixing \(\omega \in \Omega \), we define the random walk as a Markov chain - homogeneous in time - \(\left(X_{n}\right)_{n \in \mathbb{N}}\) on \(\mathbb{Z}\) or \(\mathbb{Z}^{d}\) with the following transition prosbabilities:
\begin{equation}
    p(x,y,\omega)= \mathbb{P}^{\omega}\left(X_{1}=y \vert X_{0}=0\right).
\end{equation}
\subsection{The Quenched Law and Annealed Law}
The probability \(\mathbb{P}^{\omega}\) describes the mechanisms of spatial motion of the random walk in a certain environment \(\omega\). To indicate that the random walk initially starts at point \(x\), we introduce the notation \(\mathbb{P}^{\omega}_{x}\). Thus, we have:
\begin{equation}
    p(x,y,\omega) = \mathbb{P}^{\omega}_{x}(X_{1} = y).
\end{equation}
Since the quenched law depends on the initially randomly selected environment, it defines a kind of random variable on the space of the environments. Therefore, we can take the average of this quenched probability \(\mathbb{P}^{\omega}\) with respect to the distribution of the environment represented by \(\mathbb{P}\). The result of this average builds a probability law called the annealed law \(\mathcal{P}\) such that:
\begin{equation}
    \mathcal{P}(A) = \int_{\Omega} \mathbb{P}^{\omega}(A) d \mathbb{P}(\omega).
\end{equation}

The event \(A\) may represent a certain property of the random walk. By the above equation, we understand that if \(A\) is verified almost surely (a.s.) in any environment with respect to the quenched law, it is thus verified a.s. with respect to the annealed law. 

In the Chernov model, a family of i.i.d. random variables is chosen. The chain \(X_{n}\) is only Markovian conditionally on the environment set, i.e., with respect to the quenched \(\mathbb{P}^{\omega}_{x}\). However, under the annealed law, the Markov property is lost for the chain \(X_{n}\).

To understand the loss of this Markov property, we first recall it roughly: A stochastic process \(X_{n}\) is Markovian with respect to a probability \(\mathbb{P}\) if:
\begin{equation}
\mathbb{P} \left( X_{n+1} = x_{n+1} \vert X_{n} = x_{n}, \ldots, X_{0} = x_{0} \right) = \mathbb{P} \left( X_{n+1} = x_{n+1} \vert X_{n} = x_{n} \right).
\end{equation}
The above equation ensures memory loss as the process evolves since "the future depends only on what has just preceded it."
Let's try to check it for the annealed law. First, we note:
\[
(*) = \mathbb{P} \left( X_{n+1}=x_{n+1} \vert X_{n}=x_{n}, \ldots, X_{0}=x_{0} \right)
\]
\[
(*) = \frac{\mathbb{P} (X_{n+1}=x_{n+1}, X_{n}=x_{n}, \ldots, X_{0}=x_{0})}{\mathbb{P}(X_{n}=x_{n}, \ldots, X_{0}=x_{0})}.
\]
We denote \((**)= \mathbb{P} (X_{n+1}=x_{n+1}, X_{n}=x_{n}, \ldots, X_{0}=x_{0})\). Then:
\[
(**) = \int_{\Omega} \mathbb{P}^{\omega}(X_{n+1}=x_{n+1}, X_{n}=x_{n}, \ldots, X_{0}=x_{0}) \, d \mathbb{P}(\omega).
\]
\[
(**)= \int_{\Omega} \mathbb{P}^{\omega} \left( X_{n+1}=x_{n+1} \vert X_{n}=x_{n}, \ldots, X_{0}=x_{0}\right) \mathbb{P}^{\omega}(X_{n}=x_{n}, \ldots, X_{0}=x_{0}) \, d \mathbb{P}(\omega).
\]
We apply the Markov property for \(\mathbb{P}^{\omega}\):
\[
(**) = \int_{\Omega} \mathbb{P}^{\omega} \left( X_{n+1}=x_{n+1} \vert X_{n}=x_{n}\right) \mathbb{P}^{\omega}(X_{n}=x_{n}, \ldots, X_{0}=x_{0}) \, d \mathbb{P}(\omega).
\]
Finally, we have:
\[
(*) = \frac{\int_{\Omega} \mathbb{P}^{\omega} \left( X_{n+1}=x_{n+1} \vert X_{n}=x_{n}\right) \mathbb{P}^{\omega}(X_{n}=x_{n}, \ldots, X_{0}=x_{0}) \, d \mathbb{P}(\omega)}{\int_{\Omega} \mathbb{P}^{\omega}(X_{n}=x_{n}, \ldots, X_{0}=x_{0}) \, d \mathbb{P}(\omega)}.
\]

We can see that the past of the random walk is taken into account at the next step because we can't simplify by \(\mathbb{P}^{\omega}(X_{n}=x_{n}, \ldots, X_{0}=x_{0})\) in the denominator and numerator. We understand in this way that random walking learns about the environment by making more jumps. The trace of this learning is represented by \(\mathbb{P}^{\omega}(X_{n}=x_{n}, \ldots, X_{0}=x_{0})\).
\subsection{The origin of the terms "quenched" and "annealed"}

These two words come from the field of metallurgy. They were first introduced by D. E. Temkin \cite{Temkin} in the early 1970s.

Quenching is the rapid cooling of a piece in water, oil, or air in order to obtain certain properties of the material. Being a type of heat treatment, quenching prevents undesirable low-temperature processes, such as phase transformations, from occurring. Avoiding phase transformations, for example, fixes the environment in some way. This is the idea of the quenched law.

Annealing is a heat treatment that modifies the physical and sometimes the chemical properties of a material in order to increase its ductility and reduce its hardness, making it easier to work with. During annealing, the atoms migrate in the crystal lattice, and the number of dislocations decreases, which leads to a change in ductility and hardness. In short, the environment of the material is being modified, and the name annealed law thus has a physical meaning.

\section{Concepts of recurrence and transience}
\subsection{Case of a random walk with a deterministic environment}
In order to intuitively understand recurrence and transience for a random walk, we consider a simple random walk that evolves in $\mathbb{Z}^{d}$ with transition probabilities independent of network sites so that, being in a certain $x \in \mathbb{Z}$, the walker has a probability $p$ to go right and $q=1-p$ to go left.

The two notions of recurrence and transience are a dichotomy of the answer to the question: How many times can we go back to the same starting point? If the answer is infinite, then the starting point is said to be recurrent. If the answer is a finite number of times, then the starting point is said to be transient.

In mathematical terms, we consider a sequence of independent random variables identically distributed $(Z_{i})_{i \in \mathbb{N}}$ with values of +1 or -1 such that:
\begin{equation}
    \mathcal{P}(Z_{i}= 1)= p  
\end{equation}
We will note $(X_{n})_{n \in \mathbb{N}}$ as the trajectory of the random walk as a function of time $n$. So the trajectory is written:
\begin{equation}
    X_{n}= \sum_{i=1}^{n} Z_{i} 
\end{equation}
The answer to the question of recurrence and transience requires an understanding of the asymptotic behavior of the walk $X_{n}$.
If we note $m=p-q$ as the average of the random variable $Z_{0}$, we have by the law of large numbers:
\begin{equation}
\lim_{n \to\infty} \frac{X_{n}}{n} = m.
\end{equation}
In the case of a non-symmetrical random walk, i.e., $m \neq 0$:
\begin{equation}
\lim_{n \to\infty} X_{n}= \infty .
\end{equation}
We, therefore, understand that for the non-symmetrical case of a random walk, the trajectory diverges to infinity as soon as its asymptotic velocity $v= \frac{X_{n}}{n} $ is not zero. The trajectory cannot return an infinite number of times to the starting point 0, so the point 0 is transient. 
Since all points are similar in terms of transition probability, all points are transient. In this case, we say that the random walk is transient.

We are now interested in the symmetrical case. The law of large numbers does not allow us to decide on the recurrence or transience of the walk.
However, by looking at the time of the first passage after 0 through the site $x$: $ T_{x}= \inf \{ n \geq 1 \mid X_{n}=x \}$ and by convention taking $\inf(\emptyset)= \infty$, we can show that $\mathbb{P}(T_{0} < \infty )= 1$. The random walk, therefore, almost certainly returns in finite time to its starting point of 0. Hence, 0 is recurrent. This is a property that is not exclusively reserved for 0:
\begin{thm} \ \\
For every $x \in \mathbb{Z}$, the time of the first passage $T_{x}$ is almost infinite. Once you return to the same starting point, the random walk behaves as if it will start the random walk again from the beginning. It is a result of the Markovian property. Therefore, each site is visited an infinite number of times.
\end{thm}
\subsection{In the RWRE}
We saw in the previous section that the classification of the random walk on $\mathbb{Z}$ in terms of recurrence and transience is done by evaluating the average $m=p-q$. However, this criterion is no longer valid for an RWRE in $\mathbb{Z}$.
```latex
The American mathematician Solomon \cite{Solomon} proposed a classification criterion in 1975 based on $\rho = \ln\left(\frac{q}{p}\right)$:

\begin{thm} \ \\
We define $\rho_{0} := \frac{q_{0}}{p_{0}}$ and $\eta := \mathbb{E}(\ln \rho_{0})$.

\begin{enumerate}
    \item If $\eta < 0$, then the walk is transient and $\lim_{n \to \infty} X_{n} = +\infty$ ($\mathbb{P}_{0}$-almost surely).
    
    \item If $\eta > 0$, then the walk is also transient and $\lim_{n \to \infty} X_{n} = -\infty$ ($\mathbb{P}_{0}$-almost surely).
    
    \item If $\eta = 0$, then the walk is recurrent. Moreover, we have ($\mathbb{P}_{0}$-almost surely):
    \[
    \limsup_{n \to \infty} X_{n} = +\infty \quad \text{and} \quad \liminf_{n \to \infty} X_{n} = -\infty.
    \]
\end{enumerate}
\end{thm}

\subsection{Asymptotic velocity}
When the walk is transient, it is reasonable to ask how fast it escapes to infinity. Solomon provided an answer to this question in 1975 with the following theorem \cite{Solomon}:

\begin{thm} \ \\
The limit $v := \lim_{n \to \infty} \frac{X_{n}}{n}$ exists ($\mathbb{P}_{0}$-almost surely). It is given by:
\begin{enumerate}
    \item If $\mathbb{E}(\rho_{0}) < 1$, then $v = \frac{1 - \mathbb{E}(\rho_{0})}{1 + \mathbb{E}(\rho_{0})}$.
    
    \item If $\mathbb{E}(\rho_{0}^{-1}) < 1$, then $v = -\frac{1 - \mathbb{E}(\rho_{0}^{-1})}{1 + \mathbb{E}(\rho_{0}^{-1})}$.
    
    \item Otherwise, $v = 0$.
\end{enumerate}
\end{thm}

Since the logarithm is concave, we have by Jensen's inequality: $\mathbb{E}(\ln(\rho_{0})) \leq \ln(\mathbb{E}(\rho_{0}))$. Hence, it is possible to have $\lim_{n \to \infty} X_{n} = \infty$ and $\lim_{n \to \infty} \frac{X_{n}}{n} = 0$. This result is quite surprising compared to ordinary random walks, indicating a slowdown in the transient case.

Indeed, assuming $\mathbb{E}(\rho_{0}) < 1$ and using Jensen's inequality:
\begin{equation}
    \mathbb{E}(\rho_{0}) = \mathbb{E}(\rho_{0}^{-1})^{-1} \geq \mathbb{E}(\rho_{0})^{-1} - 1,
\end{equation}
we obtain:
\begin{equation}
    0 < v \leq 2\mathbb{E}(\rho_{0}) - 1 = \mathbb{E}(p_{0} - q_{0}).
\end{equation}

We observe that the asymptotic velocity $v$ is less than $m_{0} = \mathbb{E}(p_{0} - q_{0})$, which is the asymptotic velocity in the ordinary random walk case. What is even more surprising is the possibility of zero asymptotic speed for the walk, where the trajectory advances much more slowly relative to the time scale $n$. This situation occurs when $m_{0} > 0$ and simultaneously $\eta := \mathbb{E}(\ln \rho_{0}) > 0$.

\section{ Random conductance model }
In this section, we will discuss the random conductance model that will allow us to give a probabilistic interpretation of the current and voltage in a resistance network and an intuitive understanding of the notion of recurrence/transience. This will be based on the work of Peter Doyle and Laurie Snell \cite{electric}.

\subsection{Definition of the model}
The conductance model is simply a connected graph $\mathbb{G}$ with a resistor $R_{xy}$ on each edge $xy$. Accordingly, conductivity can be defined as $C_{xy} = \frac{1}{R_{xy}}$.

We define a random walk on the graph $\mathbb{G}$ by a Markov chain with transition matrix $\mathbb{P}$ such that:
\[
\mathbb{P}_{xy} = \frac{C_{xy}}{\sum_{z} C_{xz}}.
\]

\begin{figure}[!htbp]
    \centering
    \includegraphics{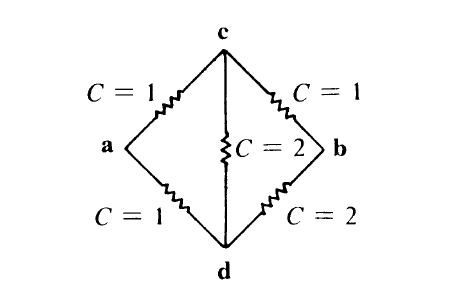}
    \caption{A conductance network}
    \label{conductance_network}
\end{figure}

\vspace{5mm}
Thus, the walker can move directly from one node to any other with probability $\mathbb{P}_{xy}$. This differs from previous models, where the walker only moves to the nearest neighboring nodes. An example of a conductance network is illustrated in Figure \ref{conductance_network}.

In the case of Figure \ref{conductance_network}, the transition matrix is given by:
\[
\mathbb{P} = \begin{bmatrix}
0 & 0 & 0.5 & 0.5 \\
0 & 0 & \frac{1}{3} & \frac{2}{3} \\
\frac{1}{4} & \frac{1}{4} & 0 & \frac{1}{2} \\
\frac{1}{5} & \frac{2}{5} & \frac{2}{5} & 0 
\end{bmatrix}.
\]
\subsection{A probabilistic interpretation of current and voltage}

Returning to our general model, since the graph is connected, the walker can move between any pair of states. A Markov chain with this property is known as an ergodic Markov chain.

Moreover, besides being ergodic, Markov chains associated with networks possess another important property: reversibility. Conversely, any ergodic and reversible Markov chain allows us to define a conductance network.

Reversibility characterizes ergodic chains in electrical networks due to the physical laws governing the behavior of direct electric currents being invariant in time.

The potential \( v_{x} \) can be interpreted as the probability that the walker returns to \( a \) before reaching \( b \), starting from \( x \).

If we initiate \( i_{a} = 1 \), the current \( i_{xy} \) represents the average net number of transitions from \( x \) to \( y \) made by the walker before returning to \( a \). These currents are proportional to the current resulting from applying a voltage of 1V between \( a \) and \( b \).

When a voltage \( v \) is applied between points \( a \) and \( b \), such that \( v_{a} = v \) and \( v_{b} = 0 \), a current \( i_{a} = \sum_{x} i_{ax} \) flows into the circuit from an external source. The amount of current flowing depends on the overall resistance of the circuit. The effective resistance \( R_{eff} \) between \( a \) and \( b \) is defined as \( R_{eff} = \frac{v_{a}}{i_{a}} \). The reciprocal quantity \( C_{eff} = \frac{1}{R_{eff}} \) is the effective conductance. If the voltage between \( a \) and \( b \) is scaled by a constant, the currents are also scaled by the same constant, indicating that \( R_{eff} \) depends only on the ratio \( \frac{v_{a}}{i_{a}} \).

We can interpret the effective conductance in a probabilistic manner as an escape probability \( p_{esc} \), which is the probability that the walker reaches \( b \) before returning to \( a \), starting from \( a \). When \( v_{a} = 1 \), the effective conductance \( C_{eff} \) equals the total current \( i_{a} \) entering \( a \). This current is given by \( i_{a} = C_{a} p_{esc} \).

Therefore:
\[
C_{eff} = C_{a} p_{esc}.
\]
\section{Recurrence as seen by the conductance model}

\subsection{Pólya Theorem}

We denote the probability that the point never returns to its starting point as \( p_{esc} \). This new \( p_{esc} \) is simply the probability of reaching infinity before returning to the starting point, with \( b = \infty \). The chain is recurrent if \( p_{esc} = 0 \), and transient if \( p_{esc} > 0 \).

In 1929, George Pólya solved the recurrence question for the simple random walk. Specifically, the answer depends only on the dimension.

\begin{thm}
The simple random walk is recurrent for dimensions \( d = 1, 2 \) and transient for dimensions \( d \geq 3 \).
\end{thm}

Here, we'll use the random conductance model to give an electrical intuition to the notion of recurrence.

\subsection{The escape probability}

We can determine the type of recurrence for an infinite network based on the properties of increasingly larger finite graphs contained within it.

Let \( G(r) \) be the graph obtained from \( \mathbb{Z}^d \) by discarding nodes whose Euclidean distance from the origin is greater than \( r \).

For norm 1 in \( \mathbb{R}^d \), we denote \( \delta G(r) \) as the "sphere" of radius \( r \) around the origin. In two dimensions, \( \delta G(r) \) is a square.

We define a random walk on \( G(r) \) as follows: the walk starts at 0 and behaves as it does on \( \mathbb{Z}^d \) until it reaches a point on \( \delta G(r) \), at which it stays. This defines a Markov chain with each point of \( \delta G(r) \) as an absorbing state. Let \( p(r)_{esc} \) be the probability that a random walk on \( G(r) \), starting at 0, reaches \( \delta G(r) \) before returning to 0. Naturally, \( p(r)_{esc} \) decreases as \( r \) increases.

The escape probability for the infinite network is then defined as:
\[
p_{esc} = \lim_{r \to \infty} p(r)_{esc}.
\]

If \( p_{esc} = 0 \), the random walk is recurrent. Otherwise, the walk is transient.

To determine \( p(r)_{esc} \) in an electrical context, we connect all points of \( \delta G(r) \) to ground and measure the current \( i(r) \) flowing in the circuit. We know that \( p(r)_{esc} = \frac{i(r)}{2d} \). If a voltage of 1 is applied, \( i(r) \) corresponds to the effective conductance between 0 and \( \delta G(r) \). In other words, \( i(r) = \frac{1}{R_{eff}(r)} \).

The effective resistance \( R_{eff} \) from the origin to "infinity" is defined as:
\[
R_{eff} = \lim_{r \to \infty} R_{eff}(r).
\]

Thus, the escape probability can be expressed as:
\[
p_{esc} = \frac{1}{2d R_{eff}}.
\]

Therefore, the random walk is recurrent if and only if the resistance from 0 to infinity is infinite. This intuitively makes sense because, by Ohm's law, the current passing through an electrical network for a fixed voltage applied at the terminals decreases as resistance increases.
\section{Excursions, traps and critical exponent}
\subsection{The excursions}
Let $T_{11}^{G}$ denote the time spent on a left excursion from site 1, which is the time required to return to 1 by making the first jump to the left.

If $\eta = \mathbb{E}(\ln \rho_{0}) < 0$, then by Solomon's theorem, $\lim_{n \to \infty} X_{n} = +\infty$, indicating that the random walk will almost surely pass through site 1 again ($\mathbb{P}_{0}$ p.s).

We define $\omega_{1} = \mathbb{E}_{1}^{\omega}(T_{11}^{G})$, representing the average ("quenched") duration of the excursion $T_{11}^{G}$. It holds that $\omega_{1} = 1 + \mathbb{E}_{0}^{\omega}(\tau_{1})$, where $\tau_{1}$ is the time needed to return to 1 after jumping to 0.

Taking the average over the environment, the annealed average duration of a left excursion is given by:
\begin{prop}
If $\mathbb{E}(\rho_{0}) < 1$, then $\mathbb{E}(\omega_{1}) = \frac{2}{1 - \mathbb{E}(\rho_{0})}$. Otherwise, $\mathbb{E}(\omega_{1}) = \infty$.
\end{prop}

Using Solomon's theorem on asymptotic velocity, where $v$ is strictly negative implying divergence towards $-\infty$, the annealed average return time to 0 is infinite.

We intuitively link $\omega_{1}$ to $\omega_{0}$. When jumping from 1 to 0, an excursion starting from 0 also moves left. Thus, returning to 1 is closely related to excursions to 0.

Considering doing $k$ excursions before jumping to 1 and avoiding further negative sites is akin to the first success after $k$ failures, described by a geometric law with probability $q_{0}^{k} p_{0}$. Given the average excursion time at 0 is $\omega_{0}$, the average excursion time at 0 before jumping to 1 is:
\[
\sum_{k=1}^{\infty} q_{0}^{k} p_{0} (k \omega_{0}).
\]

Adding the time to go from 1 to 0 at the start of the 1's excursion plus the time to go from 0 to 1 at the end of the excursion yields:
\begin{equation}
\omega_{1} = 2 + \sum_{k=1}^{\infty} q_{0}^{k} p_{0} (k \omega_{0}) = 2 + \rho_{0} \omega_{0}.
\end{equation}

The environment's translation invariance ensures $\omega_{0}$ and $\omega_{1}$ share the same distribution. By analogy, $\omega_{0} = 2 + \rho_{0} \omega_{-1}$. This implies $\omega_{0}$ depends solely on $p_{x}$ for $x \leq -1$.
\subsection{Critical exponent}

In this section, we utilize results from the article by Kesten, Kozlov, and Spitzer published in 1975 \cite{Kesten}. Chernov and Temkin observed the surprising phenomenon where a random walk can grow indefinitely with zero asymptotic velocity. This scenario arises when $\mathbb{E}(\ln(\rho_{0})) < 0$ and $\mathbb{E}(\rho_{0}) > 1$.

A natural question to ask is: What is the typical spatial order in which the random walk moves as $n$ becomes arbitrarily large?

Kolmogorov and Spitzer conjectured that $n^{-\frac{1}{k}} T_{n}$ could converge to a stable limit law, where $T_{n}$ is the first passage time to site $n$, and $k$ is the unique positive number such that:
\begin{equation}
    \mathbb{E}(\rho_{0}^{k}) = 1.
\end{equation}
This number $k$ is referred to as the critical exponent.

Kesten's paper establishes this conjecture under the condition that $\ln(\rho_{0})$ has a non-arithmetic distribution, meaning it is not concentrated on points of the form: $...,-2c,-c,0,+c,+c,+2c,...$. Additionally, it avoids the possibility of having an atom at $-\infty$.

The theorem provides convergence results in probability depending on the value of the critical exponent $k$. Specifically, for $k < 1$, the theorem states:
\begin{equation}
    \lim_{n \to \infty} \mathbb{P}\left(n^{-k} X_{n} \leq x\right) = 1 - L_{k}(x^{-\frac{1}{k}}),
\end{equation}
where $L_{k}$ is a slowly varying function. This indicates that the random walk exhibits a limit law under scaling of order $n^{\frac{1}{k}}$.

The slowdown in the random walk can be heuristically explained by the presence of "traps" in the environment. A trap can be envisioned as a region in space that is easy to enter but difficult to exit. In one dimension, this can occur with two consecutive series of sites where the transition probabilities $p_{x}$ are significantly larger on the left and smaller on the right.
\section{Sinai's localization}

\subsection{Back to the simple random walk}

We adopt the notations from Section 4.1 and return to a simple random walk on $\mathbb{Z}$. By the central limit theorem, we have the following convergence in law:
\begin{equation}
    \frac{X_{n}}{\sqrt{n}} \longrightarrow \mathcal{N}(0,1)
\end{equation}

Informally, this suggests that the average distance of the walk $X_{n}$ from the origin grows in the order of $\sqrt{n}$ as $n$ increases. This behavior is characteristic of diffusion.

In the next paragraph, we explore a surprising effect where the environment, once considered random, significantly slows down the walk compared to the ordinary random walk.

\subsection{Sinai's Theorem}

In 1982, Sinai proved a remarkable theorem \cite{Sinai} demonstrating a strong slowdown in Random Walks in Random Environment (RWRE) in the recurrent case, where $\nu = 0$.

The theorem assumes the environment is uniformly elliptical, meaning for every $x \in \mathbb{Z}$ and $\mathbb{P}$-almost surely, we have:
\begin{equation}
    0 < \delta \leq p_{x} \leq 1 - \delta < 1
\end{equation}
This condition ensures good control over the fluctuations of the environment.

\begin{thm}
Suppose the environment is uniformly elliptical and consider the recurrent case. Assuming $\sigma = \mathbb{E}(\ln^{2}(\rho_{0})) > 0$ for the environment to be truly random, we have the following convergence in law:
\begin{equation}
    \frac{\sigma X_{n}}{\ln^{2}(n)} \longrightarrow G
\end{equation}
where $G$ is a random variable.
\end{thm}

This theorem reveals that the characteristic distance of the random walk at time $n$ scales in the order of $\ln^{2}(n)$, which is significantly slower compared to the ordinary random walk. RWRE, observed on the scale of $\ln^{2}(n)$, exhibits sub-diffusive asymptotic behavior.

Additionally, there exists another remarkable result related to the previous theorem: there exists a sequence of random variables $(W_{n})$ depending only on the environment $\omega$ such that for every $\epsilon > 0$,
\begin{equation}
    \lim_{n \to \infty} \mathbb{P}_{0} \left( \left\lvert \frac{\sigma X_{n}}{\ln^{2}(n)} - W_{n} \right\rvert > \epsilon \right) = 0.
\end{equation}

By combining these convergences, we also have in law convergence $W_{n} \longrightarrow G$. What's intriguing in this second result is that RWRE, considered on the spatial scale of $\ln^{2}(n)$, becomes localized around a "random" site (determined by the environment). This phenomenon is known as Sinai's localization.
\subsection{Heuristic interpretation of the scale \texorpdfstring{$\ln^{2}(n)$}{ln2(n)}}

In the context of one-dimensional Random Walks in Random Environment (RWRE), the walker's behavior depends significantly on the site where it is located, as demonstrated by the criterion involving \( \ln(\rho_{0}) \). This dependence suggests a potential profile where each site \( z \) has an associated slope given by \( \ln(\rho_{z}) \).

We introduce a function called the random potential \( V \), defined as follows:
\[ V_{0} = 0 \quad \text{and} \quad V_{n} - V_{n-1} = \ln(\rho_{n}) \]

For the random walk \( (X_{n})_{n} \) to reach a distance \( D \), it must surpass a potential barrier of height \( \sim \sqrt{D} \):
\begin{equation}
    \max_{0 \leq i \leq j \leq D} (V_{i} - V_{j}) \sim \sqrt{D}
\end{equation}

This implies that the time required to achieve this distance scales approximately as \( t \sim \exp(\sqrt{D}) \). Therefore, setting \( t = n \) (i.e., considering the walk's evolution over time), the random walk cannot progress beyond a distance \( D \sim (\log(n))^2 \).

The motion of an RWRE can be likened to that of a particle navigating a random potential landscape. Fluctuations in the environment create potential barriers that confine the particle for extended periods.

In the regime described by Sinai's theorem, the entire environment acts as a diffusive trap. The random walk spends a significant amount of time trapped within this potential landscape. As \( n \) increases, signifying the walk's evolution over time, its movement becomes increasingly constrained. This phenomenon is akin to an "aging" effect compared to an ordinary random walk, which progresses much faster in comparison.

\subsection{The notion of local time and its link with the potential}

\begin{figure}[!htbp]
    \centering
    \includegraphics{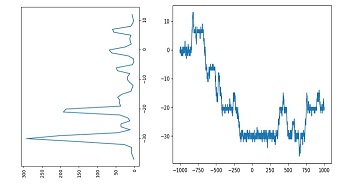}
    \caption{Simulation of the random walk (right) and its local time (left).}
    \label{MAMAetPot}
\end{figure}

The notion of local time allows us to quantify how often a random walk passes through a specific site within a fixed time interval (cf. Fig. \ref{MAMAetPot}).

For a given time \( n \), let \( \xi(x,n) \) denote the number of times the walk passes through site \( x \) between \( t=0 \) and \( t=n \).

In an ordinary random walk, the typical distance covered by the walk by time \( n \) is on the order of \( \sqrt{n} \). Chung and Hunt (1949) showed that \( \xi(x,n) \) also behaves similarly to \( \sqrt{n} \).

However, RWRE exhibits larger fluctuations in local time when compared to an ordinary random walk, especially when the typical distance is small. These fluctuations are quantified by Révész's theorem in 1986 \cite{Revesz}:

\begin{thm}
Assuming the environment is elliptical, for any \( \epsilon > 0 \), we have almost surely:
\begin{equation}
    \mathbb{P}^{\omega}_{0} \left( \xi(k,n) \leq \exp\left( \frac{\log(n)}{(\log\log(n))^{1-\epsilon}} \right) \right) \xrightarrow{n \to \infty} 1
\end{equation}

Moreover, there exists \( n_{0} < \infty \) almost surely such that for \( n \geq n_{0} \):
\begin{equation}
    \xi(k,n) \geq \exp\left( \frac{\log(n)}{(\log\log(n))^{1+\epsilon}} \right)
\end{equation}
\end{thm}

\begin{figure}[!htbp]
    \centering
    \includegraphics{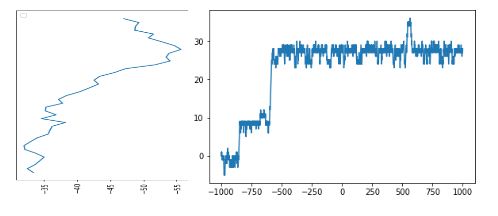}
    \caption{Simulation of random walk (right) and random potential (left).}
    \label{pot_aleat_Et_Mama}
\end{figure}

It is natural to inquire about the relationship between local time and random potential. Intuitively, deeper potential barriers imply longer durations for the walk to cross them. Figure \ref{pot_aleat_Et_Mama} illustrates a simulation of a random walk with its corresponding random potential.

As observed in Figure \ref{pot_aleat_Et_Mama}, deeper potential wells, such as the one at site 28 reaching a value of -55, cause the random walk to become highly localized around this site. Furthermore, significant variations in potential lead to the random walk "jumping" from one trap to another.

Mathematically justifying these observations, Deheuvels and Révész in 1986 \cite{Revesz} established a link between local time and the random potential:

Let \( \rho_{1} \) be the first return time of the walk to its starting point after time 0, and \( (W_{k})_{k \in \mathbb{N}} \) be a sequence depending on the random potential as follows:
\[ W_{k} = 1 + \exp(V_{2}) + \cdots + \exp(V_{k-1}) \]

The link is expressed as:
\begin{equation}
    \mathbb{P}^{\omega}_{0}(\xi(k,\rho_{1}) \geq 1) = \frac{p_{0}}{W_{k}}
\end{equation}

Equation (23) evaluates the probability that the walk passes through site \( k \) at least once before returning to its starting point. Thus, as the potential decreases towards site \( k \), \( W_{k} \) decreases, increasing the likelihood that the walk will pass through \( k \) and remain in its vicinity for some time.

\section{A conductance model applied to graphene}
\subsection{Graphene and its conduction properties}

\begin{figure}[!htbp]
    \centering
    \includegraphics{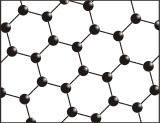}
    \caption{Graphene structure}
    \label{molecule}
\end{figure}

Graphene, a monolayer of graphite with a hexagonal structure (Cf. Fig \ref{molecule}), is a semiconductor with no energy gaps. The carrier density in graphene-based field-effect transistors (GraFETs) can be continuously adjusted from type P to type N. However, it's observed that the charge density in GraFETs exhibits mesoscopic inhomogeneity, likely due to fluctuations in trapped charges on the substrate or its surface.

For a charge-neutral graphene sheet ($ n_{g} = 0 $), this material can be envisioned as a checkerboard of doped regions (type N and P) separated by non-conductive P-N junctions. This description aligns with a random resistance model, following the findings of Cheianov et al. \cite{graphene}.

\begin{figure}[!htbp]
    \centering
    \includegraphics{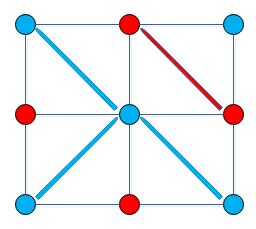}
    \caption{An example of conductance links between puddles of electrons and holes.}
    \label{P_N_Jonction}
\end{figure}

The system can be thought of as comprising "puddles" of electrons (N) and holes (P), with each puddle containing a large number of carriers. In this model, each site corresponds to either an electron puddle (marked in red) or a hole puddle (marked in blue) (Fig. \ref{P_N_Jonction}). The conductivity observed is governed by the random connections between these puddles rather than by the local conductivity of individual puddles.

Mathematically, each site \( (i,j) \) in the considered network is connected to its four neighboring sites with conductances defined as follows:

\begin{equation}
    C_{(i,j)}^{(i+1,j+1)} = \frac{c[1+(-1)^{i+j}X_{i,j}]}{2}
\end{equation}

\begin{equation}
    C_{(i,j+1)}^{(i+1,j)} = \frac{c[1-(-1)^{i+j}X_{i,j}]}{2}
\end{equation}

\begin{equation}
    C_{(i,j)}^{(i+1,j)} = C_{(i,j)}^{(i,j+1)} = \gamma c
\end{equation}

Here, \( X_{i,j} \) are random variables taking values \( \pm 1 \), each with an average \( \langle X_{i,j} \rangle = p \). The parameter \( p \) characterizes the doping level: \( p = 0 \) indicates equal probability for P-P and N-N junctions, while \( p > 0 \) (or \( p < 0 \)) indicates electron (or hole) doping. The parameter \( \gamma \ll 1 \) represents the leakage parameter, signifying that P-N and N-P junctions are less conductive.

The conductivity of the graphene sheet can be described as:

\begin{equation}
    \sigma(\gamma, p) := \langle C(L \to \infty) \rangle
\end{equation}

This expression denotes the average conductance as the size of the lattice \( L \) tends to infinity.
\subsection{Scale analysis}

\begin{figure}[!htbp]
    \centering
    \includegraphics[scale = 0.7]{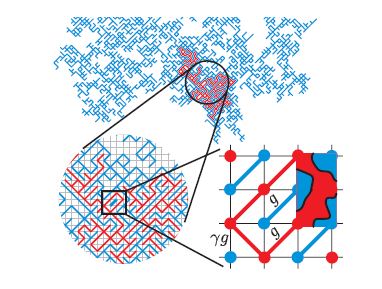}
    \caption{A layer of graphene represented by a random resistance model}
    \label{RRN_graphene}
\end{figure}

Consider a conductance model with \( p = 0 \) and \( \gamma = 0 \). The two-tone random array exhibits distinct geometric characteristics reminiscent of percolation theory. The model is critical, implying self-similar geometry at all length scales: an \( L \times L \) array typically contains large clusters of one polarity (e.g., the red clusters in the central part of Fig. \ref{RRN_graphene}) separating pairs of smaller blue clusters. Each of these clusters further contains several smaller red clusters, demonstrating scale invariance. Consequently, as \( L \) increases, increasingly larger portions of the network become fragmented due to monochromatic percolation, leading to a decrease in the observed conductance \( C(L) \). Simulation results suggest:

\begin{equation}
    \langle C(L) \rangle \sim \left(\frac{a}{L}\right)^x
\end{equation}

where \( x = 0.97 \).

\begin{figure}[!htbp]
    \centering
    \includegraphics{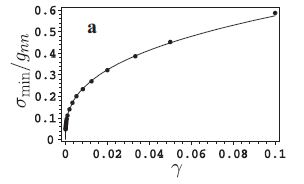}
    \caption{Evolution of conductivity depending on the leakage parameter}
    \label{sigmamiiiiin}
\end{figure}

To describe the conductivity behavior, realizations of the conductance model were generated by varying \( p \) in the interval \( [-1, 1] \) with a step size of \( \Delta p = 0.0125 \). The conductivity is then modeled as:

\begin{equation}
    \sigma(\gamma, p) = u c \gamma^\alpha F\left(\frac{p}{p_0}\right)
\end{equation}

where \( F \) is a function determined by smoothing the simulation data, and \( u \), \( \alpha \), and \( x \) are coefficients derived from these simulations. An example of the conductivity behavior obtained from simulations is shown in Fig. \ref{sigmamiiiiin} where \( \sigma_{\text{min}}(p) := \sigma(\gamma, p=0) \).

The points represent simulation results, while the curve represents the function described by Eq. (28) fitted to smooth these simulation points, illustrating the conductivity variation with respect to the leakage parameter \( \gamma \).

\section{Conclusion}

The study of Random Walks in Random Environments (RWRE) reveals how environmental variability profoundly impacts stochastic processes. Solomon's theorem underscores recurrence properties essential for understanding RWRE dynamics. Critical exponent theory by Kesten, Kozlov, and Spitzer elucidates sub-diffusive behaviors, notably in the Sinai regime, where RWRE exhibits localized behavior around specific sites.

Introducing local time illustrates the "aging" phenomenon, where walks linger near potential barriers. Applying these insights to a graphene conductance model shows how environmental fluctuations influence conductivity, offering critical perspectives on mesoscopic systems.

In conclusion, RWRE provides a robust framework for exploring randomness and environmental effects in complex systems, with implications spanning materials science and statistical physics.

%\begin{thebibliography}{9}

\newpage

\bibliography{refs}

\bibliographystyle{plain} 

% We choose the "plain" reference style 

%\end{thebibliography}
\end{document}